\documentclass[preprints,article,accept,oneauthor,pdftex,10pt,a4paper]{mdpi} 

\firstpage{1} 
\pubvolume{xx}
\issuenum{1}
\articlenumber{5}
\pubyear{2018}
\copyrightyear{2018}
\history{}%{Received: date; Accepted: date; Published: date}
\usepackage{mathtools}
\mathtoolsset{centercolon,showonlyrefs,showmanualtags}
\usepackage{dsfont}
\usepackage{nicefrac}
\DeclarePairedDelimiter{\abs}{\lvert}{\rvert}

\DeclarePairedDelimiter{\bra}{(}{)}
\DeclarePairedDelimiter{\pra}{[}{]}
\DeclarePairedDelimiter{\set}{\{}{\}}

\DeclareMathAlphabet{\mathup}{OT1}{\familydefault}{m}{n}
\newcommand{\dx}[1]{\mathop{}\!\mathup{d} #1}
\newcommand{\pderiv}[3][]{\frac{\mathop{}\!\mathup{d}^{#1} #2}{\mathop{}\!\mathup{d} #3^{#1}}}
\newcommand{\R}{\mathds{R}}
\newcommand{\cH}{\ensuremath{\mathcal H}}

\newcommand{\cN}{\ensuremath{\mathcal N}}

\DeclareMathOperator{\supp}{supp}
\DeclareMathOperator{\PI}{PI}
\DeclareMathOperator{\LSI}{LSI}

\DeclareMathOperator{\cov}{Cov}
\DeclareMathOperator{\var}{Var}
\DeclareMathOperator{\Ent}{Ent}
\DeclareMathOperator{\Id}{Id}
\DeclareMathOperator{\Diam}{diam}
\DeclareMathOperator{\osc}{osc}
\DeclareMathOperator{\Mean}{E}
\newcommand{\ttCM}{\ensuremath{\mathtt{CM}}}

\newcommand{\dsOne}{\ensuremath {\mathds{1}}}

\Title{Poincaré and log-Sobolev inequalities for mixtures}

\Author{André Schlichting$^{1,\dagger}$\orcidA{}}

\AuthorNames{André Schlichting}

\address{%
$^{1}$ \quad RWTH Aachen, Institut für Geometrie und Praktische Mathematik; schlichting@igpm.rwth-aachen.de }

\corres{Correspondence: schlichting@igpm.rwth-aachen.de }

\abstract{This work studies mixtures of probability measures on $\R^n$ and gives bounds on the Poincaré and the log-Sobolev constant of two-component mixtures provided that each component satisfies the functional inequality, and both components are close in the $\chi^2$-distance. The estimation of those constants for a mixture can be far more subtle than it is for its parts. Even mixing Gaussian measures may produce a measure with a Hamiltonian potential possessing multiple wells leading to metastability and large constants in Sobolev type inequalities. In particular, the Poincaré constant stays bounded in the mixture parameter whereas the log-Sobolev may blow up as the mixture ratio goes to $0$ or $1$. This observation generalizes the one by Chafa\"i and Malrieu to the multidimensional case. The behavior is shown for a class of examples to be not only a mere artifact of the method.}

\keyword{Poincaré inequality; log-Sobolev inequality; relative entropy; Fisher information; Dirichlet form; mixture; finite Gaussian mixtures.}

\begin{document}
\section{Introduction}

A mixture of two probability measures $\mu_0$ and $\mu_1$ on $\R^n$ is for some parameter $p\in [0,1]$ the probability measure $\mu_p$ defined by
\begin{equation}\label{mix:def:mup}
 \mu_p := p \mu_0 + (1-p) \mu_1 \:. 
\end{equation}
Hereby, both measures $\mu_0$ and $\mu_1$ are assumed to be absolutely continuous with respect to the Lebesgue measure and their supports are nested, i.e. $\supp \mu_0 \subseteq \supp \mu_1$ or $\supp \mu_1 \subseteq \supp \mu_0$. Under these assumptions at least one measure is absolutely continuous to the other one
\begin{equation}
\mu_0 \ll \mu_1 \qquad\text{or}\qquad \mu_1 \ll \mu_0  \:,
\end{equation}
which implies that that at least one of the measures has a density with respect to the other one
\begin{equation}
 \dx{\mu_0} = \pderiv{\mu_0}{\mu_1} \dx{\mu_1} \qquad\text{or}\qquad \dx{\mu_1} = \pderiv{\mu_1}{\mu_0} \dx{\mu_0} \:.
\end{equation}
This work establishes a simple and easy to check criteria under which a mixture of measures satisfies a Poincaré or log-Sobolev inequality provided that each of the component satisfies one. 
\begin{Definition}[$\PI(\varrho)$ and $\LSI(\alpha)$]\label{def:sglsi}
A probability measure $\mu$ on $\R^n$ satisfies the \emph{Poincaré inequality} with constant $\varrho > 0$, if for all functions $f:\R^n \to \R$
 \[
  \var_\mu\pra*{f} := \int \abs[\Big]{f-\int f\dx{\mu}}^2 \dx{\mu} \leq \frac{1}{\varrho} \int \abs{\nabla f}^2 \dx{\mu} \:. \tag*{$\PI(\varrho)$} 
 \]
 A probability measure $\mu$ satisfies the \emph{log-Sobolev inequality} with constant $\alpha >0$, if for all functions $f:\R^n \to \R^+$ holds
 \[
  \Ent_\mu\pra*{f} := \int f\, \log f \dx{\mu} - \int f \dx{\mu} \; \log \bra*{\int f \dx{\mu}} \leq \frac{1}{\alpha} \int \frac{\abs{\nabla f}^2}{2 f}  \dx{\mu} \:. \tag*{$\LSI(\alpha)$}
 \]
 By the change of variable $f\mapsto f^2$ the log-Sobolev inequality $\LSI(\alpha)$ is equivalent to
 \begin{equation}\label{def:lsi:square}
  \Ent_\mu\pra*{f^2} \leq \frac{2}{\alpha} \int \abs{\nabla f}^2 \dx{\mu}  \:.
 \end{equation}
\end{Definition}
The question of how $\varrho_p$ and $\alpha_p$ in $\PI(\varrho_p)$ and $\LSI(\alpha_p)$ depend for a mixture $\mu_p$ on the parameter $p\in [0,1]$ was first studied by Chafa\"i and Malrieu~\cite{Chafai2010} for measures on $\R^n$. The aim is to deduce simple criterions under which the measure $\mu_p$~\eqref{mix:def:mup} satisfies $\PI(\varrho_p)$ and $\LSI(\alpha_p)$ knowing that $\mu_{0}$ and $\mu_1$ satisfy $\PI(\varrho_{0})$, $\PI(\varrho_{1})$ and $\LSI(\alpha_{0})$, $\LSI(\alpha_1)$, respectively. The approach by Chafa\"i and Malrieu~\cite{Chafai2010} is based on a functional depending on the distribution function of the measures $\mu_0$ and $\mu_1$, which then lead to bounds on the Poincaré and log-Sobolev constant of the mixture in one dimension.

This work generalizes part of the results from Chafa\"{\i} and Malrieu~\cite{Chafai2010} to the multidimensional case by a simple argument. The estimates on the Poincaré and log-Sobolev constant hold for the case, where the $\chi^2$-distance of $\mu_0$ and $\mu_1$ is bounded (see~\eqref{e:chi01} for its definition). For this to be true, at least one of the measures $\mu_0$ and $\mu_1$ needs to be absolutely continuous to the other, which is also a necessary condition for the mixture having connected support. The resulting bound is optimal in the scaling behavior of the mixture parameter~$p\to 0,1$, i.e.~a logarithmic blow-up behavior in $p$ for the log-Sobolev constant, whereas the Poincaré constant stays bounded. This different behavior of the Poincaré and log-Sobolev constant was also observed in the setting of metastability in~\cite[Remark 2.20]{MS14}. 

Let us first introduce the principle for the Poincaré inequality in Section~\ref{s:mix:PI} and then for the log-Sobolev inequality in Section~\ref{s:mix:LSI}. Then, the procedure is illustrated on specific examples of mixtures in Section~\ref{s:mix:ex}.

\section{Poincaré inequality}\label{s:mix:PI}

To keep the presentation concise, the following notation for the mean of a function $f:\R^n\to \R$ with respect to a measure $\mu$ is introduced
\begin{equation}\label{e:mean}
 \Mean_\mu[f] := \int f \dx\mu \:.
\end{equation}
In this way, the variance in $\PI(\varrho)$ and relative entropy in $\LSI(\alpha)$ become
\begin{equation}
 \var_\mu[f]=\Mean_\mu\pra*{\bra*{f-\Mean_\mu[f]}^2} = \Mean_\mu\pra*{f}-\bra*{\Mean_\mu[f]}^2 \quad\text{and}\quad \Ent_{\mu}[f] = \Mean_\mu\pra*{f \log f} - \Mean_\mu[f] \log\bra*{\Mean_\mu[f]} \:.
\end{equation}
Likewise, the covariance of two functions $f,g:\R^n\to \R$ is defined by
\begin{equation}
  \cov_{\mu}\pra*{ f,g} = \Mean_\mu\pra*{ \bra*{f-\Mean_\mu[f]}\bra*{g-\Mean_\mu[g]}} = \Mean_{\mu}\pra*{f \,g } - \Mean_\mu[f] \, \Mean_{\mu}[g] \:.
\end{equation}
The Cauchy-Schwarz inequality for the covariance takes now the form
\begin{equation}\label{e:CS:cov}
  \cov_\mu\pra*{f,g}\leq \var_\mu[f] \, \var_\mu[g] \:.
\end{equation}
The argument is based on an easy but powerful observation for measures~$\mu_0$ and~$\mu_1$ with joint support.
\begin{Lemma}[Mean-difference as covariance]\label{lem:sg:meandiffcov}
  If $\supp \mu_0 = \supp \mu_1$, then for any $\vartheta\in [0,1]$ and any function $f:\R^n \to \R$ holds 
 \begin{equation}\label{mix:e:rep:meandiff}
  \Mean_{\mu_0}[f] - \Mean_{\mu_1}[f] = - \vartheta \cov_{\mu_0}\pra*{f,\pderiv{\mu_1}{\mu_0}} + (1-\vartheta)\cov_{\mu_1}\pra*{f,\pderiv{\mu_0}{\mu_1}} .
 \end{equation}
\end{Lemma}
\begin{proof}
 The change of measure formula yields that the covariances above are just the difference of the expectation on the right-hand side
\[ 
\cov_{\mu_0}\pra*{f,\pderiv{\mu_1}{\mu_0}} = \Mean_{\mu_0}\pra*{f\;\pderiv{\mu_1}{\mu_0}} - \Mean_{\mu_0}[f]\, \Mean_{\mu_0}\pra*{\pderiv{\mu_1}{\mu_0}} = \Mean_{\mu_1}[f] - \Mean_{\mu_0}[f] \]
 and likewise for $\cov_{\mu_1}\pra*{f,\pderiv{\mu_1}{\mu_0}}$.
\end{proof}
% \begin{Remark}\label{mixtures:rem:1dFunctional}
%  The above Lemma was the first observation in generalizing~\cite[Lemma 4.3]{Chafai2010} to the multidimensional case, but demanding that both measures are absolutely continuous to each other. In~\cite{Chafai2010} an optimal control of the mean-difference in one-dimension in terms of the distribution functions $F_i$ of the measures $\mu_i$, not necessarily absolutely continuous to each other, was deduced:
%  \begin{equation}\label{mixtures:rem:1dFunctional:e}
%   \bra*{\Mean_{\mu_0}(f) - \Mean_{\mu_1}(f)}^2 \leq I(p) \int \abs{f'}^2 \dx{\mu_p} , \quad\text{where } I(p) := \int \frac{(F_1(x)-F(0))^2}{p \mu_1(x)+q \mu_0(x)} \dx{x} . 
%  \end{equation}
%  Let us investigate under which conditions a strategy using the representation~\eqref{mix:e:rep:meandiff} leads to good results.
% \end{Remark}
The subsequent strategy is based on~\eqref{mix:e:rep:meandiff} by using a Cauchy-Schwarz inequality to arrive at the product of two variances. Then, $\PI(\varrho_0)$ or $\PI(\varrho_1)$ can be applied and the parameter $\vartheta$ leaves freedom to optimize the resulting expression. This allows to prove the following theorem, which is the generalization of~\cite[Theorem 4.4]{Chafai2010} to the multidimensional case for the Poincaré inequality provided $\mu_0$ and $\mu_1$ are absolutely continuous to each other.
\begin{Theorem}[$\PI$ for absolutely continuous mixtures]\label{thm:sg:mixtures}
  Let $\mu_0$ and $\mu_1$ satisfy $\PI(\varrho_0)$ and $\PI(\varrho_1)$, respectively, and let both measures be absolutely continuous to each other. Then for all $p \in (0,1)$ and $q=1-p$ the mixture measure $\mu_p = p \, \mu_0 + q \, \mu_1$ satisfies $\PI(\varrho_p)$ with 
\begin{equation}\label{thm:sg:mixtures:equ}
 \frac{1}{\varrho_p}\leq \begin{dcases}
  \frac{1}{\varrho_0}\:, &\text{ if } \frac{\varrho_1}{\varrho_0}\geq 1+p \chi_{1} \\
  \frac{1}{\varrho_1}\:, &\text{ if } \frac{\varrho_0}{\varrho_1}\geq 1+q \chi_{0} \\
  \frac{p \,\chi_{1} + p\,q\, \chi_{0} \,\chi_{1}  + q\,\chi_{0}}{\varrho_0\,p \,\chi_{1}  + \varrho_1 \, q \, \chi_{0} } \:, & \text{ else} . 
\end{dcases}
\end{equation}
 where 
 \begin{equation}\label{e:chi01}
  \chi_{0} := \var_{\mu_0}\pra*{\pderiv{\mu_1}{\mu_0}} \qquad \text{and} \qquad \chi_{1} := \var_{\mu_1}\pra*{\pderiv{\mu_0}{\mu_1}} . 
 \end{equation}
\end{Theorem}
\begin{proof}
 The variance of $f$ with respect to~$\mu_p$ is decomposed to
 \begin{equation}\label{thm:sg:mixtures:proof:splitvar}
   \var_{\mu_p}\pra{f} = p \var_{\mu_0}\pra{f} + q \var_{\mu_1}\pra{f} + p \, q \bra*{\Mean_{\mu_0}[f]-\Mean_{\mu_1}[f]}^2 \:.
 \end{equation}
 Hereby, the first two terms are just the expectation of the conditional variances. The second term is the variance of a Bernoulli random variable. Now, the mean-difference is rewritten by Lemma~\ref{lem:sg:meandiffcov} and the square is estimated with the Young inequality introducing an additional parameter $\eta>0$
 \[ (a+b)^2 \leq (1+\eta)\, a^2 + (1+\eta^{-1})\, b^2 \:.  \]
 Then, the Cauchy-Schwarz inequality is applied to the covariances to obtain
 \begin{equation}\label{thm:sg:mixtures:proof:est}
\begin{split}
    \var_{\mu}[f] &\leq p \var_{\mu_0}[f] + q \var_{\mu_1}[f] +\\
    &\qquad + p\, q \bra*{(1+\eta) \, \vartheta^2 \, \cov^2_{\mu_0}\pra*{f,\pderiv{\mu_1}{\mu_0}} + \bra*{1+\eta^{-1}} \, (1-\vartheta)^2\,\cov_{\mu_1}\pra*{f,\pderiv{\mu_0}{\mu_1}}} \\
   &\leq \bra*{1+(1+\eta)\,\vartheta^2 \, q \, \chi_{0}} \, p \, \var_{\mu_0}[f] + \bra*{1+\bra*{1+\eta^{-1}}\,(1-\vartheta)^2 \, p \, \chi_{1}} \, q \, \var_{\mu_1}[f] \\
   &\leq \frac{1+(1+\eta)\, \vartheta^2 \, q\, \chi_{0}}{\varrho_0} \int \abs{\nabla f}^2 \, p\dx{\mu_0} + \frac{1+\bra*{1+\eta^{-1}}\,(1-\vartheta)^2 \, p\, \chi_{1}}{\varrho_1} \int  \abs{\nabla f}^2 \, q\dx{\mu_1}\\
   &\leq \max\set[\Bigg]{\frac{1+(1+\eta)\,\vartheta^2 \, q \, \chi_{0}}{\varrho_0},\frac{1+\bra*{1+\eta^{-1}}\,(1-\vartheta)^2 \, p \, \chi_{1}}{\varrho_1}} \int  \abs{\nabla f}^2 \dx{\mu} \:.
   \end{split}
\end{equation}
The resulting maximum is now minimized in $\eta>0$ and $\vartheta\in [0,1]$. To do so without loss of generality $\varrho_0\geq \varrho_1$ is assumed. The other case can always be obtained by interchanging the roles of $\mu_0$ and $\mu_1$. If $\varrho_0 > \varrho_1$, then $\vartheta =1$ and $\eta \to 0$ is optimal as long as
  \[ \frac{1+ q \, \chi_{0}}{\varrho_0} \leq \frac{1}{\varrho_1} \: , \]
  This corresponds to the second case in~\eqref{thm:sg:mixtures:equ}. By symmetry the first case follows if $\varrho_1\geq \varrho_0$. 

  Now, in the case $\varrho_0\geq \varrho_1$ and $\varrho_0\leq (1+q \chi_{0}) \varrho_1$ there exists by monotonicity for every $\vartheta\in (0,1)$ a unique $\eta_\ast=\eta_\ast(\vartheta)>0$ such that both terms in the $\max$ of the right-hand side in~\eqref{thm:sg:mixtures:proof:est} are equal and hence the max is minimal. Since $q \, \chi_{0}>0$ and $p \, \chi_{1} >0$, the sum of the coefficients in front is then given by $h(\vartheta):=(1+\eta)\vartheta^2+(1+\tfrac{1}{\eta})(1-\vartheta)^2$ in $\vartheta$ as a function of $\eta$. The minimization of $h$ in $\vartheta\in(0,1)$ leads to $\vartheta_\ast = \tfrac{1}{1+\eta}$ and it holds
  \[ h(\vartheta^*)=\frac{1}{1+\eta} + \frac{\eta}{1+\eta} = 1 \: . \]
  Hence, in this caes the parameter $s=(1+\eta_\ast)\, \vartheta_\ast^2 = \frac{1}{1+\eta_\ast} \in (0,1)$ and $(1+\eta_\ast^{-1})\, (1-\vartheta_\ast)^2=\frac{\eta_\ast}{1+\eta_\ast} =1-s$. Thus, the problem can be rephrased: Find $s_\ast\in (0,1)$ which solves
  \[\frac{1+s \, q\, \chi_{0}}{\varrho_0}=\frac{1+(1-s) \, p\, \chi_{1}}{\varrho_1} \: .\]
  The solution $s_\ast$ is given by
  \[ s_\ast = \frac{(1+p \, \chi_{1})\varrho_0 - \varrho_1}{ \varrho_0  \, p \, \chi_{1}  + \varrho_1 \, q\, \chi_{0} } \:. \]
  For this value of $s_\ast$ the value of the $\max$ in~\eqref{thm:sg:mixtures:proof:est} is given by
  \[ \frac{1+s_\ast \, q \, \chi_{0}}{\varrho_0} = \frac{p\, \chi_{1} + \frac{\varrho_1}{\varrho_0} \, q\, \chi_{0}+(1+p\ ,\chi_{1})\, q\, \chi_{0}-\frac{\varrho_1}{\varrho_0} \, q\, \chi_{0}}{ \varrho_0 \, p\, \chi_{1}  + \varrho_1 \,q \, \chi_{0} } =  \frac{p \, \chi_{1} + \, p \, q \, \chi_{0} \,  \chi_{1}  + q\, \chi_{0}}{ \varrho_0 \, p\, \chi_{1}  + \varrho_1 \,q \, \chi_{0} } \:. \qedhere \]  
\end{proof}
\begin{Remark}
 The constants $\chi_{0}$ and $\chi_{1}$ can be rewritten if $\mu_0$ and $\mu_1$ are mutual absolutely continuous as
 \begin{equation}
  \chi_{0} =  \int \bra*{\frac{\dx\mu_1}{\dx\mu_0}}^2 \dx{\mu_0}-1  = \int \frac{\dx\mu_1}{\dx\mu_0} \dx{\mu_1} -1 \qquad\text{and}\qquad \chi_{1}= \int \bra*{\frac{\dx\mu_0}{\dx\mu_1}}^2 \dx{\mu_1} -1 = \int \frac{\dx\mu_0}{\dx\mu_1} \dx{\mu_0} -1 \:.
 \end{equation}
 This quantity is also known as \emph{$\chi^2$-distance} on the space of probability measures (cf. \cite{Gibbs2002}). The $\chi^2$-distance is a rather weak distance and therefore bounds many other probability distances. Among them is also the relative entropy. Indeed, by the concavity of the logarithm and the Jensen inequality follows
 \begin{equation}
  \Ent_{\mu_0}\pra*{\pderiv{\mu_1}{\mu_0}} = \int \log \frac{\mu_1}{\mu_0} \dx{\mu_1} \leq \log \bra[\bigg]{\int \frac{\mu_1}{\mu_0}\dx\mu_1} = \log (1 + \chi_{0}) \leq \chi_{0} ,
 \end{equation}
\end{Remark}
\begin{Remark}\label{rem:sg:bound3rdcase}
 The proof of Theorem~\ref{thm:sg:mixtures} shows that the expression for $\tfrac{1}{\varrho}$ in the last case of~\eqref{thm:sg:mixtures:equ} can be bounded above and below by
 \begin{equation}\label{rem:sg:bound3rdcase:e}
  \max\set*{\frac{1}{\varrho_0},\frac{1}{\varrho_1}} \leq \frac{p \, \chi_{1} + p\, q\, \chi_{0}\, \chi_{1}  + q\,\chi_{0}}{\varrho_0 \, p \,\chi_{1}  + \varrho_1\, q \,  \chi_{0} } \leq \max\set*{\frac{1+q \, \chi_{0}}{\varrho_0},\frac{1+p \,\chi_{1}}{\varrho_1}} \:.
 \end{equation}
 In the case, where $\chi_{0}=\chi_{1}=\chi$, the formula for $\varrho_p$~\eqref{thm:sg:mixtures:equ} simplifies to
 \begin{equation}\label{e:sg:sym}
  \frac{1}{\varrho_p} \leq \frac{1+ p\,q \,\chi}{p \, \varrho_0 + q\, \varrho_1} \:. 
 \end{equation}
\end{Remark}

\begin{Corollary}\label{cor:sg:nested}
 Let $\mu_0 \ll \mu_1$ and $\mu_0$, $\mu_1$ satisfy $\PI(\varrho_0)$, $\PI(\varrho_1)$, respectively. Then for all $p\in [0,1]$ with $q=1-p$ the mixture measure $\mu_p = p \, \mu_0 + q\, \mu_1$ satisfies $\PI(\varrho_p)$ with 
 \begin{equation}
  \frac{1}{\varrho_p} = \max\set*{\frac{1}{\varrho_0},\frac{1+p \, \chi_{1}}{\varrho_1}} .
 \end{equation}
 Likewise, if $\mu_1 \ll \mu_0$, then it holds
 \begin{equation}
  \frac{1}{\varrho_p} = \max\set*{\frac{1}{\varrho_1},\frac{1+q\, \chi_{0}}{\varrho_0}} .
 \end{equation}
\end{Corollary}
\begin{proof}
  The proof is a simple consequence of Lemma~\ref{lem:sg:meandiffcov} with $\vartheta=0$ and a similar line of estimates as in~\eqref{thm:sg:mixtures:proof:est}.
\end{proof}

\section{Log-Sobolev inequality}\label{s:mix:LSI}

In this section a criterion for $\LSI(\alpha)$ is established. It will be convenient, to establish it in the form~\eqref{def:lsi:square}. For a function $g: \R^n \to \R^+$ and two probability measures $\mu_0$ and $\mu_1$ the averaged function $\bar g: \set{0,1}\to \R^+$ is defined by
\begin{equation}\label{e:avg:f}
  \bar g(0) := \Mean_{\mu_0}[g] \quad\text{and}\quad \bar g(1) := \Mean_{\mu_1}[g] \ . 
\end{equation}
Moreover, the mixture of two Dirac measures $\delta_0$ and $\delta_1$ is by slight abuse of notation denoted by $\delta_p := p \, \delta_0 + q\, \delta_1$ for $p\in [0,1]$ and $q=1-p$. Then, the entropy of the mixture $\mu_p = p\, \mu_0 + q\, \mu_1$ is given by 
\begin{equation}\label{e:split:entropy}
  \Ent_{\mu_p}[f^2] = p \Ent_{\mu_0}[f^2] + q \Ent_{\mu_1}[f^2] + \Ent_{\delta_p}\pra*{ \overline{f^2}} \:. 
\end{equation}
The following discrete log-Sobolev inequality for a Bernoulli random variable is used to estimate the entropy of the averaged function $\bar f$. The optimal log-Sobolev constant was found by Higuchi and Yoshida~\cite{Higuchi1995} and Diaconis and Saloff-Coste~\cite[Theorem A.2.]{Diaconis1996} at the same time. 
\begin{Lemma}[Optimal log-Sobolev inequality for Bernoulli measures]\label{thm:lsi:2point}
A Bernoulli measure $\mu_p$ on $\set{0,1}$, i.e. a mixture of two Dirac measures $\delta_p = p \, \delta_0 + q \,\delta_1$ with $p \in [0,1]$ and $q=1-p$ satisfies the discrete log-Sobolev inequality
\begin{equation}\label{equ:lsi:2point}
  \Ent_{\delta_p}[g] \leq \frac{p\, q}{\Lambda(p,q)} \bra[\big]{g(0)-g(1)}^2\qquad\text{for all } g:\set{0,1}\to \R^+ \:,
\end{equation}
where $\Lambda: \R^+ \times \R^+ \to \R^+$ is the \emph{logarithmic mean} defined by
\begin{equation}
  \Lambda(p,q) :=  \frac{p-q}{\log p -\log q}, \quad\text{for $p\ne q$} \qquad\text{and}\quad \Lambda(p,p) := \lim_{q\to p} \Lambda(p,q)=p \:.
\end{equation} 
\end{Lemma}
The above result allows to estimate the coarse-grained entropy in~\eqref{e:split:entropy}. 
\begin{Lemma}[Estimate of the coarse-grained entropy]\label{lem:est:coarse-entropy}
 Let $\overline{f^2}: \set{0,1} \to \R^+$ be given by $\overline{f^2}(i) := \Mean_{\mu_i}[f^2]$ for $i\in \set{0,1}$. Then for all $p\in [0,1]$ and $q=1-p$ holds
 \begin{equation}\label{e:lem:est:coarse-entropy}
  \Ent_{\delta_p}\pra*{\overline{f^2}} \leq \frac{p\, q}{\Lambda(p,q)}\bra*{\var_{\mu_0}[f] + \var_{\mu_1}[f] +\bra*{\Mean_{\mu_0}[f] - \Mean_{\mu_1}[f]}^2} \:. 
 \end{equation}
\end{Lemma}
\begin{proof}
 Lemma~\ref{thm:lsi:2point} applied to $\Ent_{\delta_p}(\overline{f^2})$ yields
\begin{equation}\label{e:split:ent:p1}\begin{split}
   \Ent_{\bar\mu}(\overline{f^2})) &\leq  \frac{p\,q}{\Lambda(p,q)} \bra*{\sqrt{\overline{f^2}(0)}-\sqrt{\overline{f^2}(1)}}^2 \:.
  \end{split}
\end{equation}
The square-root-mean-difference on the right-hand side of~\eqref{e:split:ent:p1} can be estimated by using the fact that the function $(a,b)\mapsto (\sqrt{a}-\sqrt{b})^2$ is jointly convex on $\R^+\times \R^+$. Indeed, by introducing the functions $f_0,f_1: \R^n \times \R^n \to \R^+$ defined by $f_0(x,y)=f(x)$ and $f_1(x,y)=f(y)$, an application of the Jensen inequality yields the estimate
\begin{equation}\label{e:est:squmeandiff}
\begin{split}
    \bra*{\sqrt{\Mean_{\mu_0}\pra*{f^2}} - \sqrt{\Mean_{\mu_1}\pra*{f^2}}}^2 &= \bra*{\sqrt{\Mean_{\mu_i\times \mu_j}\pra*{f_0^2}} - \sqrt{\Mean_ {\mu_0\times\mu_1}\pra*{f_1^2}}}^2\\
    &\leq \Mean_{\mu_0\times\mu_1}\pra*{(f_0-f_1)^2} \\
    &\leq \Mean_{\mu_0}\pra[\big]{f^2} - 2\Mean_{\mu_0}\pra*{f} \, \Mean_{\mu_1}\pra*{f} + \Mean_{\mu_1}\pra[\big]{f^2} \\
    &= \var_{\mu_0}[f] + \var_{\mu_1}[f] +\bra*{\Mean_{\mu_0}[f] - \Mean_{\mu_1}[f]}^2 \:. 
   \end{split}\end{equation}
Now, a combination~\eqref{e:split:ent:p1} and~\eqref{e:est:squmeandiff} gives~\eqref{e:lem:est:coarse-entropy}.
\end{proof}
The decomposition~\eqref{e:split:entropy} together with \eqref{e:lem:est:coarse-entropy} yields that a mixture $\mu_p = p\, \mu_0 + q\,\mu_1$ for $p\in[0,1]$ and $q=1-p$ satisfies 
 \begin{equation}\label{split:equ:entropy2} \begin{split}
     \Ent_{\mu_p}\pra[\big]{f^2} &\leq  p \Ent_{\mu_0}\pra[\big]{f^2} + q \Ent_{\mu_1}\pra[\big]{f^2} \\
    &\qquad +\frac{p\, q}{\Lambda(p,q)}\var_{\mu_0}[f] + \var_{\mu_1}[f] +\bra*{\Mean_{\mu_0}[f] - \Mean_{\mu_1}[f]}^2 \:. 
    \end{split}\end{equation}
The right-hand side of~\eqref{split:equ:entropy2} consists of quantities, which can be estimated under the assumption that $\mu_0$ and $\mu_1$ satisfy $\LSI(\alpha_0)$ and $\LSI(\alpha_1)$. The following theorem provides an extension of the result~\cite[Theorem 4.4]{Chafai2010} to the multidimensional case for the log-Sobolev inequality.
\begin{Theorem}[$\LSI$ for absolutely continuous mixtures]\label{thm:lsi:mixture}
 Let $\mu_0$ and $\mu_1$ satisfy $\LSI(\alpha_0)$ resp. $\LSI(\alpha_1)$, respectively, and let both measures be absolutely continuous to each other. Then for all $p \in (0,1)$ and $q=1-p$ the mixture measure $\mu_p = p \, \mu_0 + q \, \mu_1$ satisfies $\LSI(\alpha_p)$ with 
\begin{equation}\label{thm:lsi:mixtures:equ}
 \frac{1}{\alpha_p}=\begin{dcases}
  \frac{1+q\, \lambda_p}{\alpha_0} \:, &\text{if } \frac{\alpha_1}{\alpha_0}\geq 1+p \, \lambda_p(1+ \chi_{1}(1+q \, \lambda_p)) \:, \\
  \frac{1+p\, \lambda_p}{\alpha_1} \:, &\text{if } \frac{\alpha_0}{\alpha_1}\geq 1+q \, \lambda_p(1+ \chi_{0}(1+p \, \lambda_p)) \:, \\
  \frac{p(1+q \, \lambda_p) \chi_{1} + p\,q \, \lambda_p \, \chi_{0} \, \chi_{1}  + q(1+p\, \lambda_p)\, \chi_{0}}{\alpha_0 \, p\, \chi_{1}  + \alpha_1 \, q \, \chi_{0} }\:, & \text{ else} \:.
\end{dcases}
\end{equation}
 Hereby, $\chi_{0}$ and $\chi_{1}$ are given in~\eqref{e:chi01} and $\lambda_p$ is used for the inverse logarithmic mean 
 \begin{equation}
   \lambda_p := \frac{1}{\Lambda(p,q)} = \frac{\log p - \log q}{p-q}\:,  \qquad\text{for } p\ne \frac{1}{2}\:, \qquad\text{and}\qquad \lambda_{\nicefrac{1}{2}} = 2 \:.
 \end{equation}
\end{Theorem}
\begin{proof}
 The starting point is the splitting obtained from~\eqref{split:equ:entropy2}. The variances and mean-difference in~\eqref{split:equ:entropy2} can be estimated in the same way as in the proof~\eqref{thm:sg:mixtures:proof:est} of Theorem~\ref{thm:sg:mixtures}. Additionally, the fact~\cite{Ledoux1999} that $\LSI(\alpha)$ implies $\PI(\alpha)$ is used to derive for any $\eta>0$ and any $\vartheta\in (0,1)$
 \begin{equation}\label{e:lsi:p0}
 \begin{split} 
 \Ent_{\mu_p}[f^2] &\leq \frac{1}{\alpha_0}\bra*{1+q \, \lambda_p\,\bra[\big]{1+(1+\eta)\,\vartheta^2 \,\chi_{0}}} \int \abs{\nabla f}^2 \, p \dx{\mu_0}\\ 
   &\qquad +\frac{1}{\alpha_1}\bra*{1+p \, \lambda_p \, \bra[\big]{1+(1+\eta^{-1})\, (1-\vartheta)^2 \, \chi_{1}}} \int \abs{\nabla f}^2 \, q \dx{\mu_1}\\
   &\leq \max\set*{\frac{1+q \, \lambda_p\,\bra[\big]{1+(1+\eta)\,\vartheta^2 \,\chi_{0}}}{\alpha_0} , \frac{1+p \, \lambda_p \, \bra[\big]{1+(1+\eta^{-1})\, (1-\vartheta)^2 \, \chi_{1}}}{\alpha_1} } \int \abs{\nabla f}^2 \dx{\mu_p} . \raisetag{3\baselineskip}
\end{split}
\end{equation}
 By introducing reduced log-Sobolev constants
 \begin{equation}\label{e:reduced:alpha}
 \tilde \alpha_0 := \frac{\alpha_0}{1+q \lambda_p} \qquad\text{and}\qquad \tilde \alpha_1 := \frac{\alpha_1}{1+p \lambda_p} \:, 
 \end{equation}
 as well as defining the constants $\tilde c_0$ and $\tilde c_1$ by
 \begin{equation}\label{e:reduced:chi}
   \tilde \chi_{0} := \frac{\chi_{0} \lambda_p}{1+q \lambda_p} \qquad \text{and}\qquad \tilde \chi_{1} =\frac{\chi_{1} \lambda_p}{1+p \lambda_p} \:,   
 \end{equation}
 the bound~\eqref{e:lsi:p0} takes the form
\begin{equation}\label{mix:equ:opt:LSI}
 \Ent_{\mu_p}(f^2) \leq \max\set*{\frac{1+(1+\eta)\vartheta^2 \tilde \chi_{0}}{\tilde \alpha_0} , \frac{1+(1+\tfrac{1}{\eta})(1-\vartheta)^2 \tilde \chi_{1}}{\tilde \alpha_1} } \int \abs{\nabla f}^2 \dx{\mu_p} \:.
\end{equation}
The estimate~\eqref{mix:equ:opt:LSI} has the same structure as the estimate~\eqref{thm:sg:mixtures:proof:est}, where $\tilde \alpha_0$, $\tilde \alpha_i$ play the role of $\varrho_0$, $\varrho_1$ and $\tilde \chi_{0}$, $\tilde \chi_{1}$ the roles of $\chi_{0}$, $\chi_{1}$. Hence, the optimization procedure from the proof of Theorem~\ref{thm:sg:mixtures} applies to this case and the last step consists of translating the constants $\tilde \alpha_0$, $\tilde\alpha_1$ and $\tilde \chi_{0}$, $\tilde \chi_{1}$ back to the original ones. 
\end{proof}
\begin{Remark}\label{rem:lsi:sym}
 Let the bound for $\tfrac{1}{\alpha_p}$ in the last case of~\eqref{thm:lsi:mixtures:equ} be denoted by $\frac{1}{A_p}$. Then the proof shows that it can be bounded above and below in the same way as in~\eqref{rem:sg:bound3rdcase:e} in terms of the reduced constants~\eqref{e:reduced:alpha} and~\eqref{e:reduced:chi} 
 \begin{equation}\label{rem:lsi:bound3rdcase:e}
  \max\set*{\frac{1+q \, \lambda_p}{\alpha_0},\frac{1+p \, \lambda_p}{\alpha_1}} \leq \frac{1}{A_p} \leq \max\set*{\frac{1+q\,  \lambda_p (1+ \chi_{0})}{\alpha_0},\frac{1+p \, \lambda_p (1+\chi_{1})}{\alpha_1}} \:.
 \end{equation}
 In the case $\chi_{0}=\chi_{1}=\chi$, it holds the simplified bound
 \begin{equation}\label{e:lsi:sym}
  \frac{1}{\alpha_p}\leq \frac{1+ \lambda_p + p\, q \, \lambda_p \, \chi }{p \, \alpha_0 + q\, \alpha_1} \:. 
 \end{equation}
 The inverse logarithmic mean $\lambda_p=\frac{1}{\Lambda(p,q)}$ blows up logarithmically for $p \to \set{0,1}$. Hence, even in the case $\chi=0$, the bound~\eqref{e:lsi:sym} diverges logarithmically. This logarithmic divergence looks at first sight artificial, especially in comparison to~\eqref{e:sg:sym} showing that the Poincaré constant is bounded. However, the next section with examples shows, that this blow-ups may actually occurre. Hence, the bound in~\eqref{thm:lsi:mixtures:equ} is actually optimal on this level of generality. 
\end{Remark}
An analogue statement as Corollary~\ref{cor:sg:nested} for the Poincaré constant is obtained for the lob-Sobolev constant, where the proof, following along the same lines, is omitted. 
\begin{Corollary}\label{cor:lsi:nested}
 Let $\mu_0 \ll \mu_1$ and $\mu_0$, $\mu_1$ satisfy $\LSI(\alpha_0)$, $\LSI(\alpha_1)$, respectively. 
 Then for any $p\in (0,1)$ and $p=1-q$ the mixture measure $\mu_p = p \, \mu_0 + q \, \mu_1$ satisfies $\LSI(\alpha_p)$ with 
 \begin{equation}\label{cor:lsi:nested1}
  \frac{1}{\alpha_p}\leq \max\set*{\frac{1+q \, \lambda_p}{\alpha_0}, \frac{1+p \, \lambda_p(1+\chi_{1})}{\alpha_1}}
 \end{equation} 
 Likewise, if $\mu_1 \ll \mu_0$, then it holds
\begin{equation}\label{cor:lsi:nested2}
  \frac{1}{\alpha_p}\leq \max\set*{\frac{1+p \, \lambda_p}{\alpha_1}, \frac{1+q \, \lambda_p(1+\chi_{0})}{\alpha_0}}
 \end{equation} 
\end{Corollary}
\section{Examples}\label{s:mix:ex}

The results of Theorem~\ref{thm:sg:mixtures} and Theorem~\ref{thm:lsi:mixture} are illustrated for some specific examples and also compared to the results~\cite[Section 4.5]{Chafai2010}, which however are restricted to one-dimensional measures. 
Although the criterion of Theorem~\ref{thm:sg:mixtures} and Theorem~\ref{thm:lsi:mixture} can only give upper bounds for the multidimensional case, when at least one of the mixture component is absolutely continuous to the other, it is still possible to obtain the optimal results in terms of scaling in the mixture parameter $p\to \set{0,1}$.

\subsection{Mixture of two Gaussian measures with equal covariance matrix} 

Let us consider the mixtures of two Gaussians $\mu_0:=\cN(0,\Sigma)$ and $\mu_1:=\cN(y,\Sigma)$, for some $y\in \R^n$ and a strictly positive definite covariance matrix $\Sigma\geq \sigma \Id$ in the sense of quadratic forms for some $\sigma>0$. Then, $\mu_0$ and $\mu_1$ satisfy $\PI(\sigma^{-1})$ and $\LSI(\sigma^{-1})$ by the Bakry-Émery criterion Theorem~\ref{local:thm:BakryEmery}, i.e.  $\varrho_0=\alpha_0=\varrho_1=\alpha_1 = \sigma^{-1}$. Furthermore, the $\chi^2$-distance between $\mu_0$ and $\mu_1$ can be explicitly calculated as a Gaussian integral
\[ \begin{split}
 \chi_{0}=\chi_{1}&= \frac{1}{(2\pi)^{\frac{n}{2}}\sqrt{\det{ \Sigma}}} \int \exp\bra*{-x\cdot \Sigma^{-1}x + \tfrac{1}{2}(x-y)\cdot \Sigma^{-1}(x-y)} \dx{x} -1 \\
 &=  \exp\bra*{y \cdot \Sigma^{-1} y} \frac{1}{(2\pi)^{\frac{n}{2}}\sqrt{\det{ \Sigma}}} \int \exp\bra*{-\tfrac{1}{2} (x+y)\Sigma^{-1}(x+y)} \dx{x} -1 \leq e^{\nicefrac{\abs{y}^2}{\sigma}}-1  .
\end{split}\]
Then the bound from Theorem~\ref{thm:sg:mixtures} in the form~\eqref{e:sg:sym} yields
\begin{equation}\label{mix:2Gauss:PI}
 \frac{1}{\varrho_p} \leq \bra*{1+ p\, q\, \bra*{e^{\nicefrac{\abs{y}^2}{\sigma}}-1}}\sigma.
\end{equation}
Likewise, the log-Sobolev constant follows from Theorem~\ref{thm:lsi:mixture} in the form~\eqref{e:lsi:sym} leads to
\begin{equation}\label{mix:2Gauss:LSI}
 \frac{1}{\alpha_p} \leq \bra*{1+p \, q \, \lambda_p \bra*{ e^{\nicefrac{\abs{y}^2}{\sigma}}+1}} \sigma \:. 
\end{equation}
By noting that $p q \leq p q \lambda_p \leq \frac{1}{4}$, both constants stay uniformly bounded in $p$.
The large exponential factor in the distance $e^{\nicefrac{\abs{y}^2}{\sigma}}$ cannot be avoided on this level of generality since the mixed measure $\mu_p$ has a bimodal structure leading to metastable effects~\cite[Remark 2.20]{MS14}.

The result~\cite[Corollary 4.7]{Chafai2010} deduced the following bound for $\frac{1}{\varrho_p}$ for the mixture of two one-dimensional standard Gaussians $\sigma=1$ in~\eqref{mix:2Gauss:PI}
\begin{equation}\label{mix:2Gauss:PI:CM}
 \frac{1}{\varrho_{p}} \leq 1+ p\, q \, |y|^2\bra[\Big]{\Phi(|y|) \, e^{|y|^2} + \tfrac{|y|}{\sqrt{2\pi}} \, e^{\nicefrac{|y|^2}{2}} + \tfrac{1}{2}} \: , 
\end{equation}
where $\Phi(a) = \frac{1}{\sqrt{2\pi}} \int_{-\infty}^a e^{-\nicefrac{y^2}{2}} \dx{y}$. The elementary inequalities $e^{a^2}-1 \leq a^2 e^{a^2}$ and $\Phi(a)\geq 1+\tfrac{a}{\sqrt{2\pi}} e^{-\nicefrac{a^2}{2}}$ all $a\in \R$ show that he bound~\eqref{mix:2Gauss:PI} is better than the bound~\eqref{mix:2Gauss:PI:CM} for all parameter values $p\in [0,1]$ and $|y|\geq 0$. 

Hence, this example shows, that for mixtures with components that are absolutely continuous to each other as well as whose tail behavior is controlled in terms of the $\chi^2$-distance, Theorem~\ref{thm:sg:mixtures} and Theorem~\ref{thm:lsi:mixture} even improve the bound of~\cite{Chafai2010} and generalize it to the multidimensional case.

\subsection{Mixture of a Gaussian and sub-Gaussian measure} \label{sec:mix:Gaussian:subGaussian}
Let us consider $\mu_1 = \cN(0,\Sigma)$ where $\Sigma\geq \sigma \Id$ is strictly positive definite. In addition, let the density of $\mu_0$ with respect to $\mu_1$ be bounded pointwise by some $\kappa\geq 1$, that is the relative density satisfies $\nicefrac{\dx\mu_0}{\dx\mu_1}\leq \kappa$ almost everywhere on $\R^n$. By the Bakry-Émery criterion Theorem~\ref{local:thm:BakryEmery}, it holds $\varrho_1=\alpha_1=\frac{1}{\sigma}$. Further, an upper bound for $\chi_{1}$ is obtained by the assumption on the bound on the relative density
\[ \chi_{1} = \var_{\mu_1}\bra*{\frac{\mu_0}{\mu_1}} = \int \bra*{\frac{\mu_0}{\mu_1}}^2 \dx{\mu_1} -1 \leq \kappa^2 - 1 . \]
Provided that $\mu_0$ satisfies $\PI(\varrho_0)$, the Poincaré constant of the mixture $\mu_p = p \, \mu_0 + q\, \mu_1$ satisfies by Corollary~\ref{cor:sg:nested} the estimate
\begin{equation}\label{deg:mix:equ:sg:Gauss:subGauss}
 \frac{1}{\varrho_p} \leq \max\set*{\frac{1}{\varrho_0}, (1+p (\kappa^2-1))\sigma} \:. 
\end{equation}
Similarly, Corollary~\ref{cor:lsi:nested} provides whenever $\mu_0$ satisfies $\LSI(\alpha_1)$ the following bound for the log-Sobolev constant of the mixture measure $\mu_p$
\begin{equation}\label{deg:mix:equ:lsi:Gauss:subGauss}
 \frac{1}{\alpha_p} \leq \max\set*{\frac{1+q \lambda_p}{\alpha_0}, (1+p\lambda_p \kappa^2)\sigma} . 
\end{equation}
In this case, the logarithmic blow-up of the log-Sobolev constant cannot be rules out for $p\to \set{0, 1}$, without any further information on $\mu_0$. 

\subsection{Mixture of two centered Gaussians with different variance}
For $\mu_0=\cN(0,\Id)$ and $\mu_1=\cN(0,\sigma\Id)$, the Bakry-Émery criterion Theorem~\ref{local:thm:BakryEmery} implies $\varrho_0=\alpha_0=1$ and $\varrho_1=\alpha_1= \sigma^{-1}$. The calculation of the $\chi^2$-distance can be done using the spherical symmetry and is reduced to the one dimensional integral
\begin{equation}
 \chi_{0} = \int \frac{\dx\mu_1}{\dx\mu_0}\dx{\mu_1}-1 =\frac{\cH^{n-1}(\partial B_1)}{(2\pi)^{\frac{n}{2}} \sigma^n} \int_{\R^+} r^{n-1} e^{-\bra*{\frac{1}{\sigma}-\frac{1}{2}}r^2} \dx{r} -1 \:.
\end{equation}
Hereby, $\cH^{n-1}(S^{n-1})$ denotes the $n-1$-dimensional Hausdorff measure of the sphere $\partial B_1 = \set{x\in \R^n: \abs{x}=1}$. The integral does only exist for $\sigma<2$. In this case, it can be evaluated and simplified. The bound for the constant $\chi_{1}$ follows by duality under the substitution $\sigma \mapsto \sigma^{-1}$ and is given by
\begin{equation}\label{mix:e:gaussian:chi2}
 \chi_{0} =\begin{cases}
                                  \frac{1}{\bra*{\sigma\bra*{2-\sigma}}^{\frac{n}{2}}} -1 &, \sigma < 2\\
                  + \infty &, \sigma \geq 2
                                 \end{cases}
 \qquad\text{and}\qquad \chi_{1} =  \begin{cases}
           \frac{1}{\bra*{\sigma^{-1}(2-\sigma^{-1})}^{\frac{n}{2}}}-1 &, \sigma > \frac{1}{2} \\
       + \infty  & , \sigma \leq \frac{1}{2} 
          \end{cases} \:.
\end{equation}
If $\sigma \leq \nicefrac{1}{2}$, that is for $\chi_{1}=\infty$, the bound given in Corollary~\ref{cor:sg:nested} yields
\[ \frac{1}{\varrho_p} \leq \max\set*{\sigma, 1+q \chi_{0}} = \max\set[\bigg]{\sigma,(1-q)+\frac{q}{\bra*{\sigma\bra*{2-\sigma}}^{\frac{n}{2}}}}=p+\frac{q}{\bra*{\sigma\bra*{2-\sigma}}^{\frac{n}{2}}} \]
Similarly, if $\sigma \geq 2$, that is for $\chi_{0}=\infty$, the bound becomes
\[ \frac{1}{\varrho_p} \leq \max\set*{1, (1+p \chi_{1}){\sigma}} \leq \sigma\bra[\bigg]{q+\frac{p}{\bra*{\sigma^{-1}(2-\sigma^{-1})}^{\frac{n}{2}}}} .  \]
In the case $\frac{1}{2}<\sigma < 2$, the interpolation bound~\eqref{thm:sg:mixtures:equ} of Theorem~\ref{thm:sg:mixtures} could be applied.
% \begin{equation}\label{deg:mix:equ:gaussians:cov}
%  \frac{1}{\varrho_p} \leq 
%  \begin{cases}
%    p+\frac{q}{\bra*{\sigma\bra*{2-\sigma}}^{\frac{n}{2}}} &, \sigma \leq \frac{1}{2} \\
%    \frac{p \chi_{1} + pq \chi_{0} \chi_{1}  + q\chi_{0}}{p \varrho_0 \chi_{1}  + q \varrho_1 \chi_{0} } & , \frac{1}{2}<\sigma < 2 \\
%     \sigma\bra*{q+\frac{p}{\bra*{\sigma^{-1}(2-\sigma^{-1})}^{\frac{n}{2}}}} & , 2\leq \sigma 
%  \end{cases}
% \end{equation}
% The right-hand side of~\eqref{deg:mix:equ:gaussians:cov} is continuous and the case $\frac{1}{2}<\sigma <2 $ indeed interpolates between the cases $\sigma\leq \frac{1}{2}$ and $\sigma\geq 2$. 
However, the scaling behavior for the Poincaré constant can already be observed with the estimate~\eqref{rem:sg:bound3rdcase:e} in Remark~\ref{rem:sg:bound3rdcase}, where again thanks to the symmetry $\sigma \mapsto \frac{1}{\sigma}$ it holds
\begin{equation}\label{mix:2GaussVar:PI}
 \frac{1}{\varrho_p} \leq 
 \begin{dcases}
   p+\frac{q}{\bra*{\sigma\bra*{2-\sigma}}^{\frac{n}{2}}} \:, &\text{for } \sigma \leq 1 \:, \\
    \sigma\bra[\bigg]{q+\frac{p}{\bra*{\sigma^{-1}(2-\sigma^{-1})}^{\frac{n}{2}}}} \:, & \text{for } \sigma \geq 1 \: .
 \end{dcases}
\end{equation}
Hence, the Poincaré constant stays bounded for the full range of parameter $p\in [0,1]$ and $\sigma>0$. 

% we adapt the strategy from Section~\ref{sec:mix:Gaussian:subGaussian}. Moreover, this calculation shows that the constants $\chi_{0}$ and $\chi_{1}$ are very sensitive to the tail behavior of $\mu_0$ and $\mu_1$.
% 
% We can assume w.l.o.g. that $\sigma>1$. Then, by comparing the partition sums of~$\mu_0$ and~$\mu_1$, we observe that $\mu_0\leq \sigma^{\frac{n}{2}} \mu_1$ pointwise on the level of densities. We can now apply~\eqref{deg:mix:equ:sg:Gauss:subGauss} and find for the Poincaré constant
% \begin{equation}\label{mix:2GaussVar:PI}
%  \frac{1}{\varrho_p} \leq \max\set{1, (1+p (\sigma^n-1))\sigma} = (1+p (\sigma^n-1))\sigma . 
% \end{equation}
% Similarly, by applying~\eqref{deg:mix:equ:lsi:Gauss:subGauss}, we obtain the bound for the log-Sobolev constant
% \begin{equation}\label{deg:mix:lsi:Gauss:cov}
%  \frac{1}{\alpha_p} \leq \max\set{1+q \lambda_p, (1+p\lambda_p \sigma^n)\sigma} . 
% \end{equation}
% The bound~\eqref{deg:mix:lsi:Gauss:cov} blows up logarithmically for $p\to \set{0,1}$. In the next section, we will show that one of these blow-ups is artificial and could be ruled out. In this case, it is the blow-up for $p\to 0$. 

In the case for the log-Sobolev constant, the bound from Corollary~\ref{cor:lsi:nested} gives
 \begin{equation}\label{mix:2GaussVar:LSI}
   \frac{1}{\alpha_p} \leq
 \begin{cases}
   1+\frac{q \lambda_p}{\bra*{\sigma\bra*{2-\sigma}}^{\frac{n}{2}}} &, \sigma \leq 1\\
    \sigma\bra*{1+\frac{p\lambda_p}{\bra*{\sigma^{-1}(2-\sigma^{-1})}^{\frac{n}{2}}}} & . \sigma \geq 1
 \end{cases}
\end{equation}
The bound~\eqref{mix:2GaussVar:LSI} blows up logarithmically for $p\to \set{0,1}$ in general. However, the special case $\sigma=1$, although trivially, allows for the combined bound $\frac{1}{\alpha_p} \leq 1+\min\set{p,q}\lambda_p$, which stays bounded. This behavior can be extended to the range $\sigma\in (\frac{1}{2},2)$ thanks to~\eqref{mix:e:gaussian:chi2} and the interpolation bound of Theorem~\ref{thm:lsi:mixture}. 

The result~\eqref{mix:2GaussVar:PI} can be compared with the one of~\cite[Section 4.5.2.]{Chafai2010}, which states that for some $C>0$, all $\sigma>1$ and $p\in (0,\nicefrac{1}{2})$ it holds
\begin{equation}\label{mix:2GaussVar:PI:CM}
 \frac{1}{\varrho_{p,\ttCM}} \leq \sigma + C p^\frac{1}{\sigma-1}  \:. 
\end{equation}
In general, depending on the constant $C$ the bound~\eqref{mix:2GaussVar:PI} is better for $\sigma$ small, whereas the scaling in $\sigma$ is better for~\eqref{mix:2GaussVar:PI:CM}, namely linear instead of $\sigma^{\frac{3}{2}}$ as in~\eqref{mix:2Gauss:PI}.

\subsection{Mixture of uniform and Gaussian measure}
Let $\mu_0=\cN(0,1)$ and $\mu_1=\frac{1}{\cH^n(B_1)} \dsOne_{B_1}$ with $B_1$ the unit ball around zero. The, it holds $\varrho_0=1$  by the Bakry-Émergy criterion Theorem~\ref{local:thm:BakryEmery} and $\varrho_1\geq \frac{\pi^2}{\Diam(B_1)^2}=\frac{\pi^2}{4}$ by the result of \cite{Weinberger1960}. Furthermore, since $\mu_1\ll \mu_0$ the $\chi^2$-distance between $\mu_0$ and $\mu_1$ becomes thanks to the spherical symmetry
\begin{equation}\label{mix:ex:uniGaus:c01}
\chi_{0}+1 =\int\bra*{ \frac{\mu_1}{\mu_0}}^2 \dx{\mu_0}=  \frac{\bra*{2\pi}^{\frac{n}{2}}}{\cH^n(B_1)^2} \int_{B_1} e^{\nicefrac{|x|}{2}} \dx{x} = \frac{\bra*{2\pi}^{\frac{n}{2}}\cH^{n-1}(\partial B_1)}{\cH^n(B_1)^2} \int_0^1 r^{n-1}e^{\nicefrac{r^2}{2}}\dx{r} .
\end{equation}
The volume $\cH^n(B_1)$ and the surface area $\cH^{n-1}(\partial B_1)$ of the $n$-sphere satisfy the following relations
\begin{equation}\label{mix:ex:uniGaus:gn}
\frac{\cH^{n-1}(\partial B_1)}{\cH^n(B_1)} = n \qquad\text{and}\qquad \frac{(2\pi)^{\frac{n}{2}}}{\cH^n(B_1)} = 2^{\frac{n}{2}} \Gamma\bra*{\frac{n}{2}+1}=:g_n .
\end{equation}
The integral on the right-hand side in~\eqref{mix:ex:uniGaus:c01} can be bounded below by $\frac{1}{n}$ and above by $\frac{\sqrt{e}}{n}$, which alltogether yields
\[  g_n \leq \chi_{0}+1 \leq \sqrt{e} g_n . \]
Corollary~\ref{cor:sg:nested} implies that the Poincaré constant of the mixture $\mu_p = p\, \mu_0 + q\, \mu_1$ satisfies
\begin{equation}\label{mix:ex:uniGauss:PI}
 \frac{1}{\varrho_p} \leq \max\set*{\frac{1}{\varrho_1},1+q \chi_{0}} \leq p + q\sqrt{e}\, g_n , 
\end{equation}
where the last inequality follows from $\frac{4}{\pi^2}\leq p + q\,\sqrt{e} \, g_n$ for $n\geq 1$ and all $p\in [0,1]$. 

\medskip 

The estimate the log-Sobolev constant uses that $\alpha_0 = 1$  by the Bakry-Émergy criterion Theorem~\ref{local:thm:BakryEmery} and  $\alpha_1 \geq \frac{2}{e}$ from~\eqref{e:LSI:ball}. Then, Corollary~\ref{cor:lsi:nested} yields the the bound
\begin{equation}\label{mix:ex:uniGauss:LSI}\begin{split}
 \frac{1}{\alpha_p} \leq \max\set*{\frac{1+p \lambda_p}{\alpha_1}, \frac{1+q \lambda_p(1+\chi_{0})}{\alpha_0}} 
 \leq \max\set*{\frac{(1+p \lambda_p)e}{2}, 1+q \lambda_p \sqrt{e} g_n)}.
\end{split}
\end{equation}
There is a logarithmically blow-up of the bound for $p\to\set{0,1}$. 

The blow-up for $p\to 1$ is artificial, which can be shown by a combination Bakry-Émery criterion and the Holley-Stroock perturbation principle. To do so, the Hamiltonian of $\mu_p$ is decomposed into a convex function and some error term
\begin{align}
    H_p(x) &:= -\log \mu_p(x) = -\log\bra*{\frac{p}{(2\pi)^{\frac{n}{2}}} e^{-\frac{\abs{x}^2}{2}}+ \frac{1-p}{\cH^n(B_1)} \dsOne_{B_1(0)}(x)} \notag \\
    &= -\log\bra*{ e^{-\frac{\abs{x}^2}{2}+\frac{1}{2}} + \frac{1-p}{p} \frac{(2\pi)^{\frac{n}{2}}}{\cH^n(B_1)} \sqrt{e} \dsOne_{B_1(0)}(x)} + C_{p,n} \notag  \\
    &= \frac{\abs{x}^2-1}{2} - \psi_p(x) + \tilde C_{p,n}, \label{mix:equ:uniGauss:Hp}
   \end{align}
where
\[ \psi_p(x) := \bra*{\log\bra*{ e^{-\frac{\abs{x}^2}{2}+\frac{1}{2}} + \frac{1-p}{p} \frac{(2\pi)^{\frac{n}{2}}}{\cH^n(B_1)} \sqrt{e} }+\frac{\abs{x}^2-1}{2}}\dsOne_{B_1(0)}(x)  . \]
The function $\psi_p$ is radially monotone towards the boundary of $B_1$, which yields for $\abs{x}\to 1$ the bound
\begin{equation}\label{mix:equ:uniGauss:psi}
 0\leq \psi_p(x) \leq \log\bra*{1+ \frac{1-p}{p} \frac{(2\pi)^{\frac{n}{2}}}{\cH^n(B_1)} \sqrt{e}} \:.
\end{equation}
From~\eqref{mix:equ:uniGauss:Hp} the Hamiltonian $H_p$ is compared with the convex potential $\frac{\abs{x}^2-1}{2}$ with the bound~\eqref{mix:equ:uniGauss:psi} on the perturbation $\psi_p$. This together yields by the Bakry-Émergy criterion Theorem~\ref{local:thm:BakryEmery} and the Holley-Stroock perturbation principle Theorem~\ref{local:thm:HolleyStroock} the $\mu_p$ satisfies $\PI(\tilde \varrho_p)$ and $\LSI(\tilde \alpha_p)$ with
\begin{equation}\label{mix:ex:uniGauss:BEHS}
\frac{1}{\tilde \varrho_p}\leq \frac{1}{\tilde \alpha_p} \leq 1+ \frac{1-p}{p} \sqrt{e} \; g_n , 
\end{equation}
where $g_n$ is the same constant as in~\eqref{mix:ex:uniGaus:gn}. This bound only blows up for $p\to 0$. But the blow-up is like $\frac{1}{p}$. Furthermore, the bound on the Poincaré constant is worse than the one from~\eqref{mix:ex:uniGauss:PI}. Therefore, both approaches need to be combined. 

The combination of the bounds obtained in~\eqref{mix:ex:uniGauss:LSI} and~\eqref{mix:ex:uniGauss:BEHS} results in the improved bound
\begin{equation}\label{mix:ex:uniGauss:comb}
 \frac{1}{\alpha} \leq C_n( 1 + q \lambda_p g_n ) , \quad \text{ with $C_n$ some universal constant,}
\end{equation}
which only logarithmically blows up for $p\to 0$. 

This example shows that the Poincaré constant and log-Sobolev constant may have different scaling behavior for $p\to 0$. Indeed, \cite{Chafai2010} show for this specific mixture in the one-dimensional case that the log-Sobolev constant can be bounded below by
\[ C \abs{\log p} \leq \frac{1}{\alpha} , \]
for $p$ small enough and a constant $C$ independent of $p$. In one dimension, lower bounds are accessible via the functional introduced by Bobkov-Götze~\cite{Bobkov1999a}. Hence the bound~\eqref{mix:ex:uniGauss:comb} is optimal in the one-dimensional case, which strongly indicates also optimality for the higher dimension case in terms of scaling in the mixture ration $p$.

To conclude, the Bakry-Émery criterion in combination with the Holley-Stroock perturbation principle is effective for detecting blow-ups of the log-Sobolev constant for mixtures, but has, in general, the wrong scaling behavior in the mixing parameter $p$. On the other hand, the criterion presented in Theorem~\ref{thm:lsi:mixture} provides the right scaling of the blow-up but may give artificial blow-ups, if the components of the mixture become singular in the sense of the $\chi^2$-distance.

%%%%%%%%%%%%%%%%%%%%%%%%%%%%%%%%%%%%%%%%%%
%% optional
\appendixtitles{no} %Leave argument "no" if all appendix headings stay EMPTY (then no dot is printed after "Appendix A"). If the appendix sections contain a heading then change the argument to "yes".
\appendixsections{one} %Leave argument "multiple" if there are multiple sections. Then a counter is printed ("Appendix A"). If there is only one appendix section then change the argument to "one" and no counter is printed ("Appendix").
\appendix

\section{Bakry-Émery criterion and Holley-Stroock perturbation principle}

Two classical conditions for Poincaré and log-Sobolev inequalities are stated in this part of the appendix.  The \emph{Bakry-Émery criterion} relates the convexity of the Hamiltonian of a measure and positive curvature of the underlying space to constants for the Poincaré and log-Sobolev inequality. Although the result is classical for the case of $\R^n$, the result for general convex domain was established in \cite[Theorem 2.1]{KM2016}.
\begin{Theorem}[Bakry-Émery criterion {\cite[Proposition 3, Corollaire 2]{BE},\cite[Theorem 2.1]{KM2016}}]\label{local:thm:BakryEmery}
 Let $\Omega\subset \R^n$ be convex and let $H:\Omega\to \R$ be a Hamiltonian with Gibbs measure $\mu(\dx{x})=Z_\mu^{-1} e^{-H(x)}\dsOne_{\Omega}(x) \dx{x}$ and assume that $\nabla^2 H(x) \geq \kappa >0$ for all $x\in \supp \mu$. Then $\mu$ satisfies $\PI(\varrho)$ and $\LSI(\alpha)$ with
 \begin{equation}
  \varrho \geq \kappa \qquad\text{and}\qquad \alpha \geq \kappa \:. 
 \end{equation}
\end{Theorem}
The second condition is the \emph{Holley-Stroock perturbation principle}, which allows to show Poincaré and log-Sobolev inequalities for a very large class of measures. 
\begin{Theorem}[Holley-Stroock perturbation principle {\cite[p. 1184]{HS}}]\label{local:thm:HolleyStroock}
  Let $\Omega\subset \R^n$ and $H:\Omega\to\R$ and $\psi: \Omega\to \R^n$ be a bounded function. Let $\mu$ and $\tilde \mu$ be the the Gibbs measures with Hamiltonian $H$ and $H+\psi$, respectively
 \begin{equation}
  \mu(\dx{x}) = \frac{1}{Z_\mu} e^{-H(x)} \dsOne_\Omega(x) \dx{x} \qquad\text{and}\qquad \tilde\mu(\dx{x}) = \frac{1}{Z_{\tilde \mu}} e^{-H(x)-\psi(x)} \dsOne_\Omega(x) \dx{x} \:.
 \end{equation}
 Then, if $\mu$ satisfies $\PI(\varrho)$ and $\LSI(\alpha)$, then $\tilde\mu$ satisfy $\PI(\tilde\varrho)$ and $\LSI(\tilde\alpha)$, respectively. Hereby the constants satisfy 
 \begin{equation}\label{local:e:HolleyStroockPI-LSI}
  \tilde \varrho \geq e^{-\osc \psi} \varrho \qquad\text{and}\qquad \tilde \alpha \geq e^{-\osc \psi} \alpha \:,
 \end{equation}
 where $\osc \psi := \sup_{\Omega} \psi - \inf_\Omega \psi$.
\end{Theorem}
Proofs relying on semigroup theory of Theorem~\ref{local:thm:BakryEmery} and Theorem~\ref{local:thm:HolleyStroock} can be found in the exposition by Ledoux~\cite[Corollary 1.4, Corollary 1.6 and Lemma 1.2]{Ledoux1999}. 

\begin{Example}[Uniform measure on the ball]
    The measure $\mu_1= \frac{1}{\cH^n(B_1)} \dsOne_{B_1}$, with $B_1$ is the unit ball around zero, satisfies $\LSI(\alpha_1)$ with
    \begin{equation}\label{e:LSI:ball}
      \alpha_1 \geq \frac{2}{e} \ . 
    \end{equation}
    The proof compares the measure $\mu_1$ with a family of measures 
    \[
      \nu_{\sigma}(\dx{x})=\tfrac{1}{Z_\sigma}\exp\bra*{-\sigma\abs{x}^2+\tfrac{\sigma}{2}} \dsOne(x) \dx{x} \qquad \text{for } \sigma>0 \:.
    \]
    Then, it holds that $\nu_\sigma$ satisfies $\LSI(2\sigma)$ by the Bakry-Émery criterion Theorem~\ref{local:thm:BakryEmery}. Moreover, it holds that $\osc_{x\in B_1} \abs{-\sigma |x|^2 +\sigma/2}=\tfrac{\sigma}{2}$ and hence $\mu_1$ satisfies $\LSI(2\sigma e^{-\sigma})$ by the Holley-Stroock perturbation principle Theorem~\ref{local:thm:HolleyStroock} for all $\sigma>0$. Optimizing the expression $2\sigma e^{-\sigma}$ in $\sigma$ gives the bound~\eqref{e:LSI:ball}.
\end{Example}

%%%%%%%%%%%%%%%%%%%%%%%%%%%%%%%%%%%%%%%%%%
% \funding{}

%%%%%%%%%%%%%%%%%%%%%%%%%%%%%%%%%%%%%%%%%%
\acknowledgments{This work is part of the Ph.D.~thesis~\cite{PhD12} written under the supervision of Stephan Luckhaus at the University of Leipzig. The author thanks the Max-Planck-Institute for Mathematics in the Sciences in Leipzig for providing excellent working conditions. The author thanks Georg Menz for many discussion on mixtures and metastability.}

%%%%%%%%%%%%%%%%%%%%%%%%%%%%%%%%%%%%%%%%%%
\conflictsofinterest{The author declares no conflict of interest.}

%%%%%%%%%%%%%%%%%%%%%%%%%%%%%%%%%%%%%%%%%%
% Citations and References in Supplementary files are permitted provided that they also appear in the reference list here. 

%=====================================
% References, variant A: internal bibliography
%=====================================
\reftitle{References}
\bibliographystyle{Definitions/mdpi.bst}
\bibliography{bib}

% The following MDPI journals use author-date citation: Arts, Econometrics, Economies, Genealogy, Humanities, IJFS, JRFM, Laws, Religions, Risks, Social Sciences. For those journals, please follow the formatting guidelines on http://www.mdpi.com/authors/references
% To cite two works by the same author: \citeauthor{ref-journal-1a} (\citeyear{ref-journal-1a}, \citeyear{ref-journal-1b}). This produces: Whittaker (1967, 1975)
% To cite two works by the same author with specific pages: \citeauthor{ref-journal-3a} (\citeyear{ref-journal-3a}, p. 328; \citeyear{ref-journal-3b}, p.475). This produces: Wong (1999, p. 328; 2000, p. 475)

%=====================================
% References, variant B: external bibliography
%=====================================
%\externalbibliography{yes}
%\bibliography{your_external_BibTeX_file}

%%%%%%%%%%%%%%%%%%%%%%%%%%%%%%%%%%%%%%%%%%
%% optional
% \sampleavailability{Samples of the compounds ...... are available from the authors.}

%% for journal Sci
%\reviewreports{\\
%Reviewer 1 comments and authors’ response\\
%Reviewer 2 comments and authors’ response\\
%Reviewer 3 comments and authors’ response
%}

%%%%%%%%%%%%%%%%%%%%%%%%%%%%%%%%%%%%%%%%%%
\end{document}